\documentclass[invmat]{svjour}
\usepackage{amsfonts,amssymb,amsmath,amscd,theorem, diagrams}
\spnewtheorem{Claim}{Claim}{\it}{\rm}
\spnewtheorem{notation}{Notation}{\bf}{\it}
\spnewtheorem{Convention}{Convention}{\bf}{\it}
\makeatletter
\@addtoreset{equation}{section}
\renewcommand{\theequation}%
{\thesection.\arabic{equation}}

\begin{document}

\newarrow{Corresponds}{<}{-}{-}{-}{>}
\newcommand{\eps}{{\varepsilon}}
\newcommand{\dspace}{\lineskip=2pt\baselineskip=18pt\lineskiplimit=0pt}
\newcommand{\proofend}{$\Box$\bigskip}
\newcommand{\con}{{\operatorname{Con}}}
\newcommand{\Graph}{{\operatorname{Graph}}}
\newcommand{\Span}{{\operatorname{Span}}}
\newcommand{\lcm}{{\operatorname{lcm}}}
\newcommand{\codim}{{\operatorname{codim}}}
\newcommand{\Int}{{\operatorname{Int}}}
\newcommand{\Tor}{{\operatorname{Tor}}}
\newcommand{\Ideal}{{\operatorname{Ideal}}}
\newcommand{\Spec}{{\operatorname{Spec}}}
\newcommand{\Sing}{{\operatorname{Sing}}}
\newcommand{\conv}{{\operatorname{conv}}}
\newcommand{\res}{{\operatorname{res}}}
\newcommand{\tet}{{\theta}}
\newcommand{\CC}{{\mathbb C}}
\newcommand{\FF}{{\mathbb F}}
\newcommand{\RR}{{\mathbb R}}
\newcommand{\ZZ}{{\mathbb Z}}
\newcommand{\PP}{{\mathbb P}}
\newcommand{\KK}{{\mathbb K}}
\newcommand{\DD}{{\mathbb D}}
\newcommand{\AAA}{{\mathbb A}}
\newcommand{\NN}{{\mathbb N}}
\newcommand{\GG}{{\mathbb G}}
\newcommand{\CV}{{\cal V}}
\newcommand{\CL}{{\cal L}}
\newcommand{\CH}{{\cal H}}
\newcommand{\CN}{{\cal N}}
\newcommand{\CO}{{\cal O}}
\newcommand{\CP}{{\cal P}}
\newcommand{\CM}{{\cal M}}
\newcommand{\CI}{{\cal I}}
\newcommand{\CU}{{\cal U}}
\newcommand{\CJ}{{\cal J}}
\newcommand{\CF}{{\cal F}}
\newcommand{\CS}{{\cal S}}
\newcommand{\CG}{{\cal G}}
\newcommand{\CE}{{\cal E}}
\newcommand{\cC}{{\cal C}}
\newcommand{\Del}{{\Delta}}
\newcommand{\bet}{{\beta}}
\newcommand{\kap}{{\kappa}}
\newcommand{\del}{{\delta}}
\newcommand{\sig}{{\sigma}}
\newcommand{\alp}{{\alpha}}
\newcommand{\Sig}{{\Sigma}}
\newcommand{\Gam}{{\Gamma}}
\newcommand{\gam}{{\gamma}}
\newcommand{\Lam}{{\Lambda}}
\newcommand{\lam}{{\lambda}}
\newcommand{\nek}{{,...,}}

\title{On Severi varieties on Hirzebruch surfaces}
\titlerunning{On Severi varieties on Hirzebruch surfaces}
\author{Ilya Tyomkin\thanks{The author was partially supported by the Postdoctoral fellowship provided by the Clore Foundation.}}
\institute{Department of Mathematics, Massachusetts Institute of
Technology, Cambridge, USA. \email{tyomkin@math.mit.edu}}
\date{}
\maketitle

\begin{abstract}
In the current paper we prove that any Severi variety on a
Hirzebruch surface contains a unique component parameterizing
irreducible nodal curves of the given genus in characteristic
zero.
\end{abstract}

\section{Introduction}
\begin{Convention} Throughout this paper we work over an
algebraically closed field $\KK$ of characteristic zero, and genus
always means geometric genus.
\end{Convention}
The study of Severi varieties is one of the classical problems in
algebraic geometry. Given a smooth projective surface $\Sigma$, a
line bundle $\CL\in Pic(\Sigma)$, and an integer $g$, one defines
Severi variety $V(\Sigma, \CL, g)\subset |\CL|$ to be the closure
of the locus of nodal curves of genus $g$. Then the subvariety
$V^{irr}(\Sigma, \CL, g)\subset |\CL|$ parameterizing irreducible
curves is of special interest. Originally, these varieties were
introduced by Severi (in the plane case) in order to prove the
irreducibility of the moduli spaces of curves $\CM_g$ in
characteristic zero. In Anhang F of his famous book, {\it
Vorlesungen \"uber algebraische Geometrie}, F. Severi gave a false
proof of the irreducibility of $V^{irr}(\PP^2, \CO_{\PP^2}(d),
g)$. And it took more than sixty years, till, in 1986,  Harris
proved this result \cite{H86}.

The study of various properties of Severi varieties, in particular
their degrees continued, and several formulas were obtained:
recursive formulas of Caporaso and Harris for projective plane,
and of Vakil for Hirzebruch surfaces, and a non-recursive formula
of Mikhalkin for toric surfaces. More formulas were obtained by
Kontsevich, Ran, Ruan, Tian and others. Although the degrees of
Severi varieties have been computed for any projective toric
surface, the irreducibility problem is still open in most of the
cases.

The goal of the current paper is to give a proof of the
irreducibility of $V^{irr}(\Sigma, \CL, g)\subset |\CL|$ on
Hirzebruch surfaces. We shall mention that Shevchishin announced
the same result, but up to our knowledge his argument is
incomplete. We shall also mention that the approach presented in
this paper is different from that of Shevchishin \cite{Shev04}.

The idea of the proof of the irreducibility of Severi varieties on
Hirzebruch surfaces is as follows. Let $\Sigma_n$ be a Hirzebruch
surface, and let $\CL=\CO_{\Sigma_n}(dL_0+kF)$ be a line bundle,
where $L_0$ and $F$ denote the effective classes generating
$Pic(\Sigma_n)$, satisfying $L_0.F=1$, $F^2=0$, and $L_0^2=n$. We
use the notation $V_{g,d,k}$ for the Severi variety $V(\Sigma_n,
\CL, g)$. In these notations the first part of the proof is given
by Proposition \ref{contain}, which states the following: any
irreducible component $V\subseteq V_{d,k,g}$ contains a very
reducible nodal curve of a special type. Now, the irreducibility
of $V_{g,d,k}^{irr}$ follows from Proposition \ref{unique},
claiming that there exists a unique component containing such a
curve, whose generic point corresponds to an irreducible curve.

The proof of Proposition \ref{contain} is similar to the plane
case proof of Harris presented in \cite{HM98}. The proof of
Proposition \ref{unique} is reduced to a combinatorial statement
using monodromy-type arguments.

To finish the introduction we shall mention that there are
examples of (non-rational) surfaces admitting reducible Severi
varieties. Moreover, these Severi varieties can have components of
different dimensions. So it is unclear how to characterize the
surfaces admitting only irreducible Severi varieties.
Nevertheless, we would like to state the following conjecture
motivated by our result and Mikhalkin's work \cite{Mikh05}:

\begin{conjecture}\label{conj0}
If $\Sigma$ is a toric surface, $\CL\in Pic\Sigma$ is an effective
class, and $g\ge 0$ is a non-negative integer then the Severi
variety $V_{\Sigma,\CL,g}^{irr}$ parameterizing irreducible nodal
curves of genus $g$ in $|\CL|$, that do not contain singular
points of $\Sigma$, is either empty or irreducible.
\end{conjecture}
We shall discuss the proof of this conjecture for rational curves
in the last section of the paper, and we hope that the general
case can be approached using the rapidly developing methods of
tropical geometry combined with the approach presented in the
current paper.

Finally, we would like to mention that the case of positive
characteristic is still open even for plane curves. It seems that
the main missing ingredient is a statement characterizing Severi
varieties in terms of their dimensions (similar to Theorem
\ref{th:dimbd}). There is a tropical evident that such a statement
must exist, however we do not know any algebraic theorem of this
type. Nevertheless we suppose that Conjecture \ref{conj0} is true
in arbitrary characteristic.

\begin{acknowledgement}
I am very grateful to J. Bernstein, G.-M. Greuel, E. Shustin, and
M. Temkin for helpful discussions. \newline Parts of this work
were done while the author was a Clore postdoctoral fellow at the
Weizmann Institute of Science, a Moore instructor at MIT, and was
visiting the Max-Planck-Institut f\"ur Mathematik at Bonn. I would
like to thank these institutions for their hospitality.
\end{acknowledgement}

\section{Preliminaries}

\subsection{Deformation theory}

In this section we discuss several (basic) facts from the
deformation theory of algebraic varieties and algebraic maps. Most
of the statements, ideas, and proofs presented here can be found
in different sources (see for example
\cite{AC81,CH98.1,H86,Vak00}, and \cite{Zar82} for related
topics). However, I decided to write it down here for the
completeness of the presentation.

\subsubsection{Deformations of maps.}

Let $X$ and $Y$ be smooth algebraic varieties over an
algebraically closed field $\KK$ of characteristic zero, and let
$f:X\rTo Y$ be an algebraic map. In this section we discuss the
deformation theory of the pair $(X,f)$, namely, we fix $Y$ and
vary $X$ and $f$. Denote $\DD=\Spec\KK[\epsilon]/(\epsilon^2)$. We
recall that a first-order deformation of $(X,f)$ is a triple
\begin{itemize}
\item a flat family $\pi:\widetilde{X}\rTo \DD$,
\item a map $F:\widetilde{X}\rTo Y$,
\item an isomorphism $\alpha: (\widetilde{X}_0, F_0)\rTo (X, f)$, where
$\widetilde{X}_0=\widetilde{X}/(\epsilon)$, and
$F_0=F/(\epsilon)$.
\end{itemize}

\begin{notation}
The set of first order deformations of the pair $(X, f)$ modulo
isomorphisms is denoted $Def^1(X, f)$.
\end{notation}

\begin{proposition}\label{tan:prop} If $0\rTo T_X\rTo f^*T_Y$ is
exact, then
$$Def^1(X,f)\cong H^0(X, \CN_f),$$ where $\CN_f$ denotes the normal
sheaf to $f$, i.e. the cokernel of the map $$df:T_X\rTo f^*T_Y.$$
\end{proposition}
\begin{proof} First, we choose affine coverings $Y=\cup_{i=1}^n
Y_i$ and $X=\cup_{i=1}^n X_i$ such that $f(X_i)\subset Y_i$ for
all $i$. Now, let $\xi\in Def^1(X, f)$ be a first-order
deformation. We shall use the following well-known claim

\begin{Claim} Let $Z=\Spec A$ be a smooth affine variety over an
algebraically closed field $\mathbb{K}$, and let $Z_\epsilon=\Spec
A_\epsilon$ be an infinitesimal extension of $Z$, i.e. a pair
consisting of a flat morphism $Z_\epsilon\rTo \DD$ together with
an isomorphism $Z_\epsilon/(\epsilon)\simeq Z$. Then $Z_\epsilon$
is isomorphic to the trivial extension, namely $A_\epsilon\simeq
A\oplus \epsilon A$.
\end{Claim}
Due to the claim we can fix trivializations
\begin{equation}\label{is1}
\CO_{\widetilde{X}}(X_i)\simeq \CO_{X_i}(X_i)\oplus
\epsilon\CO_{X_i}(X_i).
\end{equation}
Then we obtain the automorphisms
$$\beta_{ij}:\CO_{X_{ij}}(X_{ij})\oplus
\epsilon\CO_{X_{ij}}(X_{ij})\rTo \CO_{X_{ij}}(X_{ij})\oplus
\epsilon\CO_{X_{ij}}(X_{ij}),$$ equal to identity modulo
$\epsilon$. Hence $${\beta_{ij}(x+\epsilon
y)=x+\epsilon(y+D_{ij}(x))},$$ where
$D_{ij}:\CO_{X_{ij}}(X_{ij})\rTo \CO_{X_{ij}}(X_{ij})$ are
derivations. The maps
$$F^*:\CO_Y(Y_i)\rTo\CO_{X_i}(X_i)\oplus \epsilon\CO_{X_i}(X_i)$$
are given by $F^*(x)=f^*(x)+\epsilon D_i(x)$, where
$$D_i:\CO_Y(Y_i)\rTo \CO_{X_i}(X_i)$$ are derivations, and the
following equality holds:
$$D_i-D_j=D_{ij}\circ f^*.$$
Hence the set $D^\xi=(D_1,...,D_n)$ defines a global section of
the sheaf $\CN_f$. It is clear that $D^\xi$  does not depend on
the choice of the trivializations in (\ref{is1}). Now, one can
easily check that the constructed correspondence provides us with
the bijection\footnote{Here one must use the fact that $0\rTo
T_X\rTo f^*T_Y$ is exact.}
$$Def^1(X, f)\ni\xi \rCorresponds D^\xi\in H^0(X, \CN_f),$$ which
in fact is an isomorphism of vector spaces. Moreover, this
bijection does not depend on the choice of coverings
$Y=\cup_{i=1}^n Y_i$ and $X=\cup_{i=1}^n X_i$. \qed
\end{proof}

\subsubsection{Families of curves on algebraic surfaces.}

Let $\Sigma$ be a smooth projective algebraic surface, and let
$\CL$ be a line bundle on $\Sigma$. Consider an irreducible
variety $V\subseteq |\CL|$ whose generic element is a reduced
curve. The goal of this section is to give a natural upper bound
on the dimension of $V$. Let
$$
\begin{diagram}
  {\cal C} &        & \rInto &        & V\times\Sigma \\
           & \rdTo  &        & \ldTo  &               \\
           &        &    V   &        &               \\
\end{diagram}
$$
be the tautological family of curves over $V$, and let
$\widetilde{\cal C}\rTo \cal C$ be its normalization. Then for
almost all $p\in V$ the fiber $\widetilde{\cal C}_p$ is the
normalization of the fiber\footnote{Here we use the fact that
$\KK$ is a field of characteristic zero.} ${\cal C}_p$.

Let us choose a generic point $0\in V$. Due to the generic
flatness theorem and Proposition \ref{tan:prop}, we then have a
natural map
$$\mu : T_0V\rTo Def^1(C, f)\cong H^0(C, \CN_{f}),$$
where $C=\widetilde{{\cal C}_0}$ and $f$ is the composition of
maps $$C=\widetilde{{\cal C}_0}\rTo^{}_{} {\cal C} \rTo
\Sigma\times \{0\}=\Sigma.$$

\begin{proposition}\label{prop:tansp} (1) The map $T_0V\rTo H^0(C, \CN_{f})$ is an
embedding.

\noindent (2) $Im(T_0V)\cap H^0(C, \CN_f^{tor})=0$.
\end{proposition}
\begin{proof} First, we choose a smooth projective irreducible curve
$D\subset \Sigma$ intersecting $f(C)$ transversally, such that
$h^0(\Sigma, \CL(-D))=0$. Let
$$\begin{diagram}
  X &\rTo^F & \Sigma \\
  \dTo & & \\
  \DD & &
\end{diagram}$$
be a first-order deformation of the pair $(C, f)$. We define the
new pair ${X_D=X\times_{\Sigma} D}$ and $F_D:X_D\rTo D$.
\begin{Claim}\label{defth:cl1}
The natural map $X_D\rTo \DD$ is flat.
\end{Claim}
So we constructed a map
$$\rho : Def^1(C, f)\rTo Def^1(f(C)\cap D, f_D),$$ where
$f_D=F_D/(\epsilon)$.
\begin{Claim}\label{defth:cl2}
$H^0(C, \CN_f^{tor})\subset Ker (\rho)$.
\end{Claim}
To finish the proof it is enough to show that
$\dim(\rho(\mu(T_0V)))=\dim(V).$ Consider the exact sequence
$$H^0(\Sigma, \CL(-D))\rTo H^0(\Sigma, \CL)\rTo H^0(D,
\CL\otimes\CO_D).$$ The first term is zero, hence the map
 $\alpha : V\rInto |\CL\otimes\CO_D|$ is an
embedding. Now
$$\begin{array}{ll}
T_{\alpha(0)}|\CL&\otimes\:\CO_D|=\\
&H^0\left(f(C)\cap D, \CO_{f(C)\cap D}(f(C)\cap
D)\right)=Def^1(f(C)\cap D, f_D)
\end{array}
$$ and $\rho\circ\mu=d\alpha(0)$, which
implies $\dim(\rho(\mu(T_0V)))=\dim(V)$. \qed
\end{proof}

\begin{proof}[of Claim \ref{defth:cl1}] First of all we
can assume that $D, \Sigma,$ and $C$ are affine. Thus $X$ and
$X_D$ are affine as well. Next we shall use the following lemma:
\begin{lemma}[\cite{Eis95}, Proposition 6.1]
Let $R$ be a commutative ring with identity, and let $M$ be an
$R-$module. If $I$ is an ideal of $R$, then $Tor^1_R(R/I, M)=0$ if
and only if the map $I\otimes_RM\rTo M$ is an injection. The
module $M$ is flat if and only if this condition is satisfied for
every finitely generated ideal $I\subset R$.
\end{lemma}
In our case $R=\KK[\epsilon]/(\epsilon^2)$ and the only
non-trivial ideal we have is $I=\KK\epsilon$. By the lemma it is
enough to show that the map
$$\CO_{X_D}\otimes_{\CO_\DD}\KK\epsilon\rTo^\gamma \CO_{X_D}$$ is an
embedding. Since $\CO_X$ is flat over $\CO_\DD$, we have an exact
sequence
$$0\rTo \CO_X\otimes_{\CO_\DD} \KK\epsilon\rTo\CO_X\rTo\CO_C\rTo 0.$$
Tensoring this sequence with $\CO_D$ over $\CO_\Sigma$ we obtain
$$0\rTo Tor^1_{\CO_\Sigma}(\CO_D,
\CO_C)\rTo\CO_{X_D}\otimes_{\CO_\DD}\KK\epsilon\rTo^\gamma
\CO_{X_D}.$$
%\rTo \CO_D\otimes_{\CO_\Sigma}\CO_C\rTo 0.$$
It remains to prove that $Tor^1_{\CO_\Sigma}(\CO_D, \CO_C)=0$. Let
$$0\rTo \CO_\Sigma\rTo^{\cdot j}\CO_\Sigma\rTo\CO_D\rTo 0$$
be a free resolution of $\CO_D$, where $j\in\CO_\Sigma$ is a
function defining $D$. Thus
$$Tor^1_{\CO_\Sigma}(\CO_D, \CO_C)=Ker\left(\CO_C\rTo^{f^*(j)}\CO_C
\right)=0,$$ since $f(C)$ intersects $D$ transversally. \qed
\end{proof}

\begin{proof}[of Claim \ref{defth:cl2}] To prove this
claim we have to construct the map $\rho$ explicitly. Consider the
fibered product diagram
$$\begin{diagram}
  f(C)\cap D & \rTo^{\, \, \, \, \, f_D} & D \\
  \dTo^i &  & \dTo_{i_D} \\
  C & \rTo_f & \Sigma
\end{diagram}
$$
The intersection $f(C)\cap D$ is transversal and its points are
smooth on both $f(C)$ and $D$. Hence we have the following exact
sequences on $f(C)\cap D$:
$$\begin{diagram}
0 & \rTo & i^*T_C    & \rTo^{i^*(df)}   & i^*f^*T_\Sigma     & \rTo & i^*\CN_f       & \rTo & 0\\
  &      &           &             &     \|             &      &                &      &  \\
0 & \rTo & f_D^*T_D  & \rTo^{df_D} & f_D^*i_D^*T_\Sigma & \rTo & f_D^*\CN_{i_D} & \rTo & 0\\
  &      &    \|     &             &                    &      &                &      &  \\
  &      & \CN_{f_D} &             &                    &      &                &      &
\end{diagram}$$
Moreover, $$i^*(df)(i^*T_C)\oplus df_D(f_D^*T_D)=i^*f^*T_\Sigma
.$$ Hence the map ${\gamma : \CN_{f_D}\rTo i^*\CN_f}$ is an
isomorphism, and $\rho$ is given by the composition
$$H^0(C, \CN_f)\rTo^{i^*} H^0(f(C)\cap D,
i^*\CN_f)\rTo^{\gamma^{-1}} H^0(f(C)\cap D, \CN_{f_D}).$$ To
finish the proof we note that $df(p)\ne 0$ for any $p\in
i(f(C)\cap D)$. Thus $i^*\CN_f^{tor}=0$, which implies $H^0(C,
\CN_f^{tor})\subset Ker (\rho)$. \qed
\end{proof}

%\begin{remark} Both claims are true for $C$ and $D$ of arbitrary
%dimensions, but the proofs are a little more complicated.
%\end{remark}

\begin{theorem}\label{th:dimbd}
Let $\Sigma$ be a smooth projective algebraic surface, $C'\subset
\Sigma$ be a smooth curve, and let $\CL$ be a line bundle on
$\Sigma$. Consider a positive dimensional irreducible variety
$V\subseteq |\CL|$ whose generic element $C_0$ is a reduced curve.
Assume that $C_0^i.K_\Sigma<-1$ for any irreducible component
$C_0^i\subseteq C_0$. Then
\begin{equation}\label{th:cond0}
\dim(V)\le -C_0.K_\Sigma+g-1.
\end{equation}
Furthermore, if the equality holds, and
\begin{equation}\label{th:cond1}
C_0^i.K_\Sigma<-3
\end{equation}
for any singular irreducible component $C_0^i$ of $C_0$, then $V$
has no fixed points, $C_0$ has only nodes as its singularities,
and $C_0$ intersects $C'$ transversally.
\end{theorem}
\begin{proof} Let $C$ be the normalization of $C_0$ and let $f:C\rTo\Sigma$ be the natural map. Then  $\dim(h^0(C,
\CN_f/\CN_f^{tor}))\le -C_0.K_\Sigma+g-1$ implies
(\ref{th:cond0}), by Proposition \ref{prop:tansp}. So it is enough
to show the analogous inequality for every irreducible component
of $C$. Thus we can assume that $C$ is irreducible.

Choose an invertible sheaf $\CF$ on $C$ such that the sequence
$$0\rTo \CN_f/\CN_f^{tor}\rTo \CF\rTo \CN_f^{tor}\rTo 0$$
is exact (the existence of such $\CF$ is completely obvious). Then
\begin{equation}\label{th:eq1}
c_1(\CF)=c_1(\CN_f)=c_1(f^*T_\Sigma)-c_1(T_C)=2g-2+c_1(f^*T_\Sigma)>2g-1.
\end{equation}
Hence
$$h^0(C,\CN_f/\CN_f^{tor})\le h^0(C,\CF)=c_1(\CF)+1-g=-C_0.K_\Sigma+g-1$$
by the Riemann-Roch theorem, and the equality holds if and only if
$\CN_f^{tor}=0$.

For the second part we note that if the dimension of $V$ equals
$-C_0.K_\Sigma+g-1$, then $\CN_f^{tor}=0$, and hence $df\ne 0$
everywhere. So it remains to prove that $C_0$ has no triple
points, and that all its double points have two different tangent
directions. If $p\in C_0$ is a triple point and $q_1, q_2, q_3\in
C$ are three points mapped to $p$, then any $\xi\in T_0V$
vanishing at $q_1$ and $q_2$, must vanish at $q_3$ as well.
However, due to the Riemann-Roch theorem, inequality
(\ref{th:eq1}), and condition (\ref{th:cond1}), there exists
$\eta\in H^0(C,\CN_f/\CN_f^{tor})$ such that
$\eta(q_1)=\eta(q_2)=0\ne\eta(q_3)$. Thus
$$\dim T_0V<h^0(C,\CN_f/\CN_f^{tor})=-C_0.K_\Sigma+g-1,$$
which contradicts the equality in (\ref{th:cond0}). If $p\in C_0$
is a double point with a unique tangent direction and $q_1, q_2\in
C$ are the two pre-images of $p$, then any $\xi\in T_0V$ vanishing
at $q_1$, must also vanish at $q_2$. However, applying
Riemann-Roch theorem, inequality (\ref{th:eq1}), and condition
(\ref{th:cond1}), we can find $\eta\in H^0(C,\CN_f/\CN_f^{tor})$,
such that $\eta(q_1)=0\ne\eta(q_2)$. Hence
$$\dim
T_0V<h^0(C,\CN_f/\CN_f^{tor})=-C_0.K_\Sigma+g-1,$$ which is a
contradiction.

It remains to prove that if $\dim(V)=-C_0.K_\Sigma+g-1$, then
$C_0$ intersects $C'$ transversally. Assume that
$\dim(V)=-C_0.K_\Sigma+g-1$. Then the system $V$ has no fixed
components; hence $C_0$ does not contain $C'$. If $p$ is iether a
point of a non-transversal intersection of $C_0\cap C'$ or a fixed
point of $V$, then either $p$ has at least two pre-images
$q_1,q_2\in C$ or any $ \xi\in T_0V$ vanishes at $q$, where $q\in
C$ is the unique pre-image of $p$. In the first case any $\xi\in
T_0V$ vanishing at $q_1$ must also vanish at $q_2$. So both cases
contradict the Riemann-Roch theorem, due to (\ref{th:eq1}) and
(\ref{th:cond1}). \qed
\end{proof}

\begin{lemma}\label{lem:dimbd}
Let $(\Sigma, \CL)$ be a smooth rational surface equipped with a
line bundle, and let $L\subset \Sigma$ be a smooth curve. Let
$p_1,...,p_r\in L$ be arbitrary points, $k_1,...,k_r$ be
non-negative integers, and let
$$R\subset\left\{C\in |\CL|\, :\,  L\cap
C=\sum_{i=1}^rk_ip_i\right\}$$ be a non-empty subvariety. Choose a
generic curve $C\in R$. If $C$ is reduced and
$C_i.(K_\Sigma+L)<-1$ for any irreducible component $C_i\subset
C$,  then
$$\dim(R)\le -C.K_\Sigma+g(C)-1-L.C.$$ Moreover, if the equality holds,
and for any singular irreducible component $C_i$ of $C$ we have
$C_i.(K_\Sigma+L)<-3$, then $R$ has no fixed points but
$p_1,...,p_r$, $C$ has only nodes as its singularities outside of
$L$, for any smooth irreducible curve $C'$, not tangent to $L$, a
generic curve $D\in R$ intersects $C'$ transversally at any
$q\notin\{p_1,...,p_r\}$, and for any $i$, $D$ is a union of
smooth branches in a neighborhood of $p_i$, not tangent to $C'$.
\end{lemma}
\begin{proof} The proof is by induction on $L.C$. If
$L.C=0$, then the lemma follows from Theorem \ref{th:dimbd}, since
$L.C_i=0$ for any irreducible component $C_i\subseteq C$.

Assume now that $k_i>0$  for some $i$. Without loss of generality
$k_1>0$. Consider the blow up
$\widetilde{\Sigma}=Bl_{p_1}(\Sigma)$ with its natural projection
$\pi:\widetilde{\Sigma}\rTo \Sigma$. We denote the strict
transform of $L$ by $\widetilde{L}$. Then
$\CO_{\widetilde{\Sigma}}(\widetilde{L})\simeq\pi^*\CO_\Sigma(L)\otimes
\CO_{\widetilde{\Sigma}}(-E)$, where $E$ denotes the exceptional
divisor $\pi^{-1}(p_1)$. Let $\widetilde{C}\subset
\widetilde{\Sigma}$ be the strict transform of $C$. Then
$\CO_{\widetilde{\Sigma}}(\widetilde{C})\simeq\pi^*\CO_\Sigma(C)\otimes
\CO_{\widetilde{\Sigma}}(-m_1E)$, where $m_1=mult_{p_1}(C)$.
Define $p'_1=E\cap \widetilde{L}$ and $p'_i=p_i $ for any $i>1$,
and consider the pullback of (an open dense subset of) $R$
$$R_1\subset\left\{X\in |\CO_{\widetilde{\Sigma}}(\widetilde{C})|\, :\, g(X)=g(C)\, \mbox{and}\, X\cap \widetilde{L}=\sum_{i=1}^rk'_ip'_i\right\},$$
where $k'_1=m_{p'_1}(\widetilde{C},\widetilde{L})$ and $k'_i=k_i$
for all $i>1$. Since
$K_{\widetilde{\Sigma}}\equiv\pi^*K_\Sigma+E$,
$$\widetilde{C}_i.(K_{\widetilde{\Sigma}}+\widetilde{L})=C_i.(K_\Sigma+L)<-1$$
for any irreducible component $\widetilde{C}_i\subset
\widetilde{C}$. Since $\widetilde{L}.\widetilde{C}=L.C-m_1<L.C$,
$$\begin{array}{lll}
\dim(R)&=&\dim(R_1)\\
&\le&
-\widetilde{C}.K_{\widetilde{\Sigma}}+g(\widetilde{C})-1-\widetilde{L}.\widetilde{C}=-C.K_\Sigma+g(C)-1-L.C,
\end{array}$$
by the induction hypothesis. Moreover, if the equality holds and
for any singular component $C_i$ of $C$ we have
$C_i.(K_\Sigma+L)<-3$, which implies
$\widetilde{C}_i.(K_{\widetilde{\Sigma}}+\widetilde{L})<-3$, then
$\widetilde{C}$ has only nodes as its singularities outside of
$\widetilde{L}$, for any smooth irreducible curve $C'$ not tangent
to $\widetilde{L}$, a generic curve $D_1\in R_1$ intersects $C'$
transversally at any $q'\notin\{p'_1,...,p'_r\}$, and for any $i$,
$D_1$ is a union of smooth branches in a neighborhood of $p'_i$
not tangent to $C'$. Thus no germ is tangent to the exceptional
divisor $E$, hence $D=f(D_1)$ satisfies the required properties,
since $R_1$ has no fixed points but $p_1',...,p_r'$. \qed
\end{proof}

\subsection{Severi varieties on Hirzebruch surfaces}

Let $\Sigma_n={\PP}roj(\CO_{\PP^1}\oplus\CO_{\PP^1}(n))$ be the
Hirzebruch surface and let $\pi:\Sigma_n\rTo \PP^1$ be the natural
projection. Consider two sections $(1,0), (0, \sigma)\in
H^0(\PP^1, \CO_{\PP^1}\oplus\CO_{\PP^1}(n))$. They define the maps
$$\left(\PP^1\backslash Z(\sigma)\right)\rTo \Sigma_n.$$ We denote
the closures of the images of these maps by $L_0$ and $L_\infty$,
respectively. It is clear that $L_\infty$ is independent of the
choice of $\sigma$. The following facts will be useful
\begin{itemize}
\item The Picard group $Pic(\Sigma_n)$ is a free abelian group generated by the classes $F$ and $L_\infty$,
where $F$ denotes the fiber of the projection $\pi$. It is
important to mention that $L_0\equiv nF+L_\infty$.
\item The intersection form on $NS(\Sigma_n)=Pic(\Sigma_n)$ is given by $F^2=0$, $L_\infty^2=-n$, and $F.L_\infty=1$.
\item Any effective divisor $M\in Div(\Sigma_n)$ is linearly
equivalent to a linear combination of $F$ and $L_\infty$ with
non-negative coefficients. Moreover, if $M$ does not contain
$L_\infty$, then it is linearly equivalent to a combination of $F$
and $L_0$ with non-negative coefficients.
\item The canonical class is
$$K_{\Sigma_n}\equiv -(2L_\infty+(2+n)F)\equiv
-(L_0+L_\infty+2F).$$
\item Any smooth curve $C\equiv dL_0+kF$ has genus
$g(C)=\frac{(d-1)(dn+2k-2)}{2}.$
\end{itemize}
Now let us define the Severi varieties on $\Sigma_n$.
\begin{definition} (1) Let
$$\breve{\Sigma}_n^\delta=
\left\{(p_1,...,p_\delta)\in\Sigma_n^\delta\, \mid\, p_i\ne p_j \,
\mbox{for any}\,\, i\ne j\right\}$$ be the configuration space of
$\delta$ points in $\Sigma_n$. For non-negative integers $d, k,
\delta$, we define the decorated Severi variety
$$U_{d,k,\delta}\subset |\CO_{\Sigma_n}(dL_0+kF)|\times
\breve{\Sigma}_n^\delta$$ to be
$$\left\{(C; p_1, ... , p_\delta)\, \mid\,
C\, \mbox{is reduced,}\, L_\infty\nsubseteq C,\, \mbox{and}\,\,
p_1,...,p_\delta\in C\, \mbox{are nodes}\right\}.$$

\noindent (2) Let $g, d, k$ be non-negative integers. We define
the Severi variety $V_{g,d,k}\subseteq |\CO_{\Sigma_n}(dL_0+kF)|$
to be the closure of the locus of reduced nodal curves of genus
$g$ which do not contain $L_\infty$, and we define
$V_{g,d,k}^{irr}\subset V_{g,d,k}$ to be the union of the
irreducible components whose generic points correspond to
irreducible curves.
\end{definition}
Next, we establish the basic properties of (decorated) Severi
varieties:
\begin{proposition}\label{prop:smoothness} (1) Let $d,k,\delta$ be non-negative integers.
Then either $U_{d,k,\delta}$ is empty or $U_{d,k,\delta}$ is a smooth
equidimensional variety of dimension
$$\frac{(d+1)(nd+2k+2)}{2}-1-\delta=\dim|\CO_{\Sigma_n}(dL_0+kF)|-\delta.$$
(2) Let $\psi : U_{d,k,\delta}\rTo |\CO_{\Sigma_n}(dL_0+kF)|$ be
the projection to the first factor. Then for
$g=\frac{(d-1)(nd+2k-2)}{2}-\delta$
$$V_{g,d,k}=\left\{%
\begin{array}{ll}
    \overline{\psi\left(U_{d,k,\delta}\right)} & \hbox{if \, $\delta\ge 0$;} \\
    \emptyset & \hbox{otherwise.} \\
\end{array}%
\right.$$
\end{proposition}
\begin{proof} (1) Let us choose arbitrary $(C_0; p_{10}, ... , p_{\delta 0})\in
U_{d,k,\delta}$. We can find an open subset $U\subset \Sigma_n$
isomorphic to $\AAA^2$ and containing all the points $p_{10}, ...
, p_{\delta0}$. Fix a trivialization of
$\CO_{\Sigma_n}(dL_0+kF)(U)\simeq \KK[x,y]$. Then in a
neighborhood of $(C_0; p_{10}, ... , p_{\delta 0})$
$U_{d,k,\delta}$ is given by the (homogeneous in the coefficients)
equations $f(p_l)=f'_x(p_l)=f'_y(p_l)=0$, $1\le l\le\delta$, where
$f\in H^0(\Sigma_n, \CO_{\Sigma_n}(dL_0+kF))$. We denote $(x_{l0},
y_{l0})=p_{l0}$ and define $a_{lij}$ to be the coefficients of
$f(x+x_{l0}, y+y_{l0})$, $f=\sum a_{ij}x^iy^j \in \KK[x,y]$. Let
$f_0(x,y)=\sum \beta_{ij}x^iy^j$ be an equation defining $C_0$.
Then
$$d(f(p_l))_{|_{(f_0, p_{10}, ... , p_{\delta 0})}}=da_{l00},$$
$$d(f'_x(p_l))_{|_{(f_0, p_{10}, ... , p_{\delta 0})}}=da_{l10}+2\beta_{l20}dx_l+\beta_{l11}dy_l,$$ and
$$d(f'_y(p_l))_{|_{(f_0, p_{10}, ... , p_{\delta 0})}}=da_{l01}+\beta_{l11}dx_l+2\beta_{l02}dy_l.$$
The points $p_{10}, ... , p_{\delta 0}$ are nodes of $C_0$, hence
the matrices
$$\begin{pmatrix}
  2\beta_{l20} & \beta_{l11} \\
  \beta_{l11} & 2\beta_{l02}
\end{pmatrix}$$
are invertible. So, it remains to prove that
$da_{100}=...=da_{\delta 00}=0$ defines a subspace of codimension
$\delta$ in the tangent space to $|\CO_{\Sigma_n}(dL_0+kF)|$ at
$C_0$. In other words we have to prove that
\begin{equation}\label{eqsevcl}
h^0(C_0, \CJ(dL_0+kF))=h^0(C_0, \CO_{C_0}(dL_0+kF))-\delta,
\end{equation}
where $\CJ$ denotes the ideal of the zero-dimensional reduced
subscheme $$X=\cup_{i=1}^\delta p_{i0}\subset C_0.$$ Thus the
following Claim implies (\ref{eqsevcl}).

\begin{Claim}\label{cl:h1van1}
$H^1(C_0, \CJ(C_0))=0$.
\end{Claim}
We postpone the proof of the claim till the end of the proof of
the proposition.

(2) The inclusion $$V_{g,d,k}\subseteq\left\{%
\begin{array}{ll}
    \overline{\psi\left(U_{d,k,\delta}\right)} & \hbox{if \, $\delta\ge 0$;} \\
    \emptyset & \hbox{otherwise.} \\
\end{array}%
\right.$$ is obvious. Let us prove the opposite direction. If
$g\ge\frac{(d-1)(nd+2k-2)}{2}$, then there is nothing to prove.
Otherwise let $U\subseteq U_{d,k,\delta}$ be an irreducible
component, and let $(C;p_1,...,p_\delta)\in U$ be a generic point.
Then by part (1) and Theorem \ref{th:dimbd}
\begin{equation}\label{cl:ineq}
-C.K_{\Sigma_n}+g(C)-1\ge\dim\psi\left(U_{d,k,\delta}\right)=\dim|\CO_{\Sigma_n}(dL_0+kF)|-\delta.
\end{equation}
Hence $g(C)\ge \frac{(d-1)(nd+2k-2)}{2}-\delta=g_{ar}(C)-\delta$,
which is possible only if the equality holds. Thus the equality
holds in (\ref{cl:ineq}) as well, and by Theorem \ref{th:dimbd}
this implies the nodality of the curve $C$. So
$$\psi(C;p_1,...,p_\delta)\in V_{g,d,k}$$ and we are done.
\qed
\end{proof}

\begin{proof}[of Claim \ref{cl:h1van1}]
Let us denote the irreducible components of $X=C_0$ by $X_i$, and
the normalizations of $X_i$ by $\widetilde{X}_i$. Then the
normalization $\widetilde{X}$ of $X$ is the disjoint union of
$\widetilde{X}_i$. Consider the conductor ideal
$\CJ^{cond}=Ann(\oplus\CO_{\widetilde{X}_i}/\CO_X)\subseteq
\oplus\CO_{\widetilde{X}_i}=\CO_{\widetilde{X}}$ and the direct
sum of the conductor ideals $\oplus \CJ^{cond}_i=\oplus
Ann(\CO_{\widetilde{X}_i}/\CO_{X_i})\subseteq
\CO_{\widetilde{X}}$. Let us notice that $\CJ^{cond}\subseteq
\CO_X\subseteq\CO_{\widetilde{X}}$ is an ideal in both algebras
$\CO_X$ and $\CO_{\widetilde{X}}$, hence it is sufficient to prove
that $H^1(X, \CJ^{cond}(X))=H^1(\widetilde{X}, \CJ^{cond}(X))=0$,
since
$$0\rTo \CJ^{cond}(X)\rTo \CJ(X)_{|_{X}}$$
is exact, and the factor is a torsion sheaf.

Consider the exact sequence
$$0\rTo \CJ^{cond}\rTo \oplus \CJ^{cond}_i \rTo \oplus\CF_i\rTo 0,$$
where $\CF_i$ are torsion sheaves supported at the preimages of
the points of intersection with other irreducible components of
$X$. Let us estimate the degree of $\CF_i$. To do this we can
assume that the surface and the curves are affine. Let $f_i=0$ be
an equation of $X_i$. Then $\prod_{j\ne
i}f_j\cdot\CJ^{cond}_i\subseteq \CJ^{cond}$. Thus
$$\deg(\CF_i)\le (X-X_i).X_i=X.X_i-X_i^2.$$

It is well known that $\deg(\CJ^{cond}_i)\ge -2\delta(X_i)$ where
$\delta(X_i)$ denotes the total delta invariant of $X_i$, moreover
the equality holds for singular curves on smooth surfaces . Thus
$$\begin{array}{ll}\deg(\CJ^{cond}(X)_{|_{\widetilde{X}_i}})=&\deg(\CJ_i^{cond}(X))-\deg(\CF_i)\ge\\&
X.X_i-2\delta(X_i)-(X-X_i).X_i=\\&(K_{\Sigma_n}+X_i).X_i-2\delta(X_i)-K_{\Sigma_n}.X_i=\\
&2g_{ar}(X_i)-2-2\delta(X_i)-K_{\Sigma_n}.X_i>2g(X_i)-2,\end{array}$$
by the adjunction formula, since $-K_{\Sigma_n}.X_i>0$. Applying
Riemann-Roch theorem we conclude that $H^1(\widetilde{X}_i,
\CJ^{cond}(X)_{|_{\widetilde{X}_i}})=0$ for all $i$, which implies
$H^1(\widetilde{X}, \CJ^{cond}(X))=0$. \qed
\end{proof}

It is easy to see that if $(C; p_1, ... , p_\delta)\in
U_{d,k,\delta}$ and $C$ has exactly $\delta$ nodes, then the map
$\psi : U_{d,k,\delta}\rTo V_{g,d,k}$ for
$g=\frac{(d-1)(nd+2k-2)}{2}-\delta$, is \'etale in a neighborhood
of $(C; p_1, ... , p_\delta)$.
\begin{corollary}\label{cor:tan} If $V_{g,d,k}\ne \emptyset$ then
it has pure dimension $nd+2k+2d+g-1$, and for any $C\in V_{g,d,k}$
having exactly $\delta=\frac{(d-1)(nd+2k-2)}{2}-g$ nodes,
$V_{g,d,k}$ is smooth at $C$, and $T_C(V_{g,d,k})\simeq H^0(C,
\CJ^{cond}(C))$, where $\CJ^{cond}\subset \CO_C$ is the conductor
ideal.
\end{corollary}

\begin{proposition}\label{prop:dimbd}
Let $p_1,...,p_r\in L_0\cup L_\infty$ be fixed points, let $d>0$,
$k,k_1,...,k_r\ge 0$ be integers, and let
$$R\subset\left\{D\in |\CO_{\Sigma_n}(dL_0+kF)|\, :\,  (L_0\cup L_\infty)\cap
D=\sum_{i=1}^rk_ip_i\right\}$$ be a non empty subvariety whose
generic point $C$ corresponds to a reduced irreducible curve of
genus $g$. Then $$\dim X\le
-(K_{\Sigma_n}+L_0+L_\infty).(dL_0+kF)+g-1.$$ Moreover, if the
equality holds, then $C$ has only nodes as its singularities
outside of $L_0\cup L_\infty$, for any smooth irreducible curve
$C'$ not tangent to $L_0\cup L_\infty$, a generic curve $D\in R$
intersects $C'$ transversally outside of $\{p_1,...,p_r\}$, and
for any $i$, $D$ is a union of smooth branches in a neighborhood
of $p_i$, not tangent to $C'$.
\end{proposition}

\begin{proof}
$C.(K_{\Sigma_n}+L_0+L_\infty)=-2d<-1$, and if $C$ is singular
then $d>1$, hence $C.(K_{\Sigma_n}+L_0+L_\infty)<-3$. So we can
apply Lemma \ref{lem:dimbd} to prove the Proposition. \qed
\end{proof}

\section{The Result}\label{sec:res}

\begin{theorem}\label{irth}
Let $g, k, d$ be non-negative integers. If the variety $V_{g, d,
k}^{irr}$ is not empty then it is irreducible.
\end{theorem}
\begin{proof} Let $L_i$, $1\le i\le d$, be
generic curves in the linear system $|\CO_{\Sigma_n}(L_0)|$, and
let $F_i$, $1\le i\le k$, be generic curves in the linear system
$|\CO_{\Sigma_n}(F)|$. Define
$$\Gamma=L_1\cup ... \cup L_d\cup F_1\cup ... \cup F_k\in
\left|\CO_{\Sigma_n}(dL_0+kF)\right|.$$

\begin{proposition}\label{contain} Consider an arbitrary
component $V$ of the Severi variety $V_{g, d, k}$. Then $\Gamma\in
V$.
\end{proposition}

\begin{proposition}\label{unique}
There exists a unique component $V\subset V_{g, d, k}^{irr}$
containing $\Gamma$.
\end{proposition}
The theorem now follows. \qed
\end{proof}

\begin{remark}
One can generalize the propositions above to prove the following
more general statement: {\it Let $g, k, d, m_1,...,m_r$ be
non-negative integers satisfying $\sum m_i=k$. Consider the
varieties
$$V_{g, d, k}(m_1,...,m_r)\subset V_{g, d, k},$$ parameterizing
curves having $r$ points of tangency of orders $m_1,...,m_r$ with
$L_\infty$, and $$V_{g, d, k}^{irr}(m_1,...,m_r)=V_{g, d,
k}(m_1,...,m_r)\cap V_{g,d,k}^{irr}.$$ If $V_{g, d,
k}^{irr}(m_1,...,m_r)\ne \emptyset$ then it is irreducible.} The
generalization is pretty much straightforward, but it makes the
presentation more complicated, so we will not write it down in
this paper, but rather leave to the interested reader as an
exercise.
\end{remark}

\subsection{Proof of Proposition \ref{contain}}

$V$ is birational to a product of components of Severi varieties
whose generic points correspond to irreducible curves modulo a
finite group of symmetries, due to Claim \ref{prop:smoothness} and
Theorem \ref{th:dimbd}. Thus, without loss of generality, we can
assume that the generic point of $V$ corresponds to an irreducible
curve.

Now, to prove the proposition, it is sufficient to show that $V$
contains a nodal curve $E=L\cup E'$ where $L$ is a smooth curve of
type $L_0$ and $E'\in V_{g',d',k'}=V_{g',d-1,k}.$ If $d=1$ then
$g=0$, hence $V=V_{0,1,k}=|\CO_{\Sigma_n}(L_0+kF)|$, and there is
nothing to prove. So we can assume that $d>1$.

Let $p_1^0,...,p_{nd+k+1}^0\in L_0$ and
$p_1^\infty,...,p_k^\infty\in L_\infty$ be generic points. A
generic $C\in V$ intersects $L_0\cup L_\infty$ transversally due
to Lemma \ref{lem:dimbd}, and the locus of curves in $V$
intersecting $L_0\cup L_\infty$ transversally along $C\cap
(L_0\cup L_\infty)$ has codimension $|C\cap (L_0\cup L_\infty)|$
by Proposition \ref{prop:dimbd}. Then the locus $W$ of irreducible
curves passing through
$\{p_i^0\}_{i=1}^{nd+k}\cup\{p_j^\infty\}_{j=1}^k$ has pure
dimension
$$-(K_{\Sigma_n}+L_0+L_\infty).(dL_0+kF)+g-1;$$
and the locus of curves $V_{L_0}\subset \overline{W}$ passing
through $p_{nd+k+1}^0$, i.e. containing $L_0$ as a component, has
pure dimension
$$-(K_{\Sigma_n}+L_0+L_\infty).(dL_0+kF)+g-2.$$ Consider a map from an irreducible smooth
germ curve $$j:(T,0)\rTo (\overline{W},V_{L_0}),$$ such that
$j(0)\in V_{L_0}$ is a generic point and
$j(T^*)=j(T\setminus\{0\})\subset W$. Then for any $t\in T^*$,
$C_t$ is a nodal curve of genus g containing
$p_1^0,...,p_{nd+k}^0,p_1^\infty,...,p_k^\infty$, where $C_t$
denotes the fiber over $t$ of the corresponding flat family
$\cC\rTo T$. The central fiber $C_0$ can be presented as
$C_0=s_0L_0\cup C_0'$, where $C_0'$ is a curve that does not
contain $L_0$, and $s_0\ge 1$.
\begin{lemma}\label{lm:rednotcont}
The curve $C_0'$ does not contain $L_\infty$, and $s_0=1$.
Moreover, $C_0'$ is a nodal curve, and the points of intersection
$C_0'\cap L_0$ are smooth points of $C_0'$.
\end{lemma}
\begin{proof}
Assume that $C_0=s_0L_0\cup C_0''\cup s_\infty L_\infty$, where
$C_0''$ contains neither $L_0$ nor $L_\infty$ as its components.
Then $C_0''\equiv (d-s_0-s_\infty)L_0+(k+ns_\infty)F$. After
proceeding with an appropriate base change and replacing the
family $\cC$ by its normalization, we can consider a semistable
model $\widetilde{\cC}\rTo \cC\rTo T$ of the family $\cC$, whose
total space is smooth, generic fiber is also smooth and has genus
$g$, and its central fiber is a nodal curve. Let $\widetilde{C}_0$
be the central fiber of the semistable family, and let
$f:\widetilde{C}_0\rTo C_0$ be the natural map. Then
$\widetilde{C}_0=A_0\cup B\cup A_\infty$, where $A_i\subset
f^{-1}(L_i)$ are the unions of the connected components of
$\widetilde{C}_0$ mapped surjectively onto $L_i$, $B$ is the union
of all other components, and the following equality
holds\footnote{if $s_\infty=0$, i.e. $A_\infty=\emptyset$, then
$p_a(A_\infty):=1$.}
$$p_a(A_0)+p_a(B)+p_a(A_\infty)-2+|A_0\cap B|+|A_\infty\cap
B|=g.$$

Next step is to estimate the degrees of freedom of $C_0''$. Let
$D\subset C_0''$ be any irreducible component with reduced
structure, and let $B_D$ be any irreducible component of $B$
mapped surjectively on $D$. Then, by Proposition \ref{prop:dimbd},
$D$ vary in a family of dimension at most
$$-D.(K_{\Sigma_n}+L_0+L_\infty)+g(D)-1+l_D,$$
where $l_D$ is the number of points of $D\cap (L_0\cup L_\infty)$
distinct from $\{p_j^i\}$.

{\it Case 1:} $D\ne F$. Since
$-D.(K_{\Sigma_n}+L_0+L_\infty)=2D.F>0$, $D$ vary in a family of
dimension at most
$$-B_D.f^*(K_{\Sigma_n}+L_0+L_\infty)+g(B_D)-1+l_{B_D},$$
where $l_{B_D}$ is the number of points mapped onto $(L_0\cup
L_\infty) \setminus \{p_j^i\}$. Moreover, the equality holds if
and only if $B_D\rTo D$ is the normalization map, $l_{B_D}=l_D$,
$D$ is nodal away from $L_0\cup L_\infty$, and all its branches
are smooth at the points of intersection with $L_0\cup L_\infty$.

{\it Case 2:} $D=F$. In this case $D$ also varies in a family of
dimension at most
$$-B_D.f^*(K_{\Sigma_n}+L_0+L_\infty)+g(B_D)-1+l_{B_D}.$$
Since $p_j^i\in L_i$ are general then $1\le l_D\le 2$. If $l_D=1$
then there is a point $q\in B_D$ mapped onto one of $\{p_j^i\}$,
such that $q$ is a smooth point of $B$. Thus the pullback of $L_i$
to $B$ is reduced at this point hence $B_D\rTo D$ is an
isomorphism.

The only points of $B$ that are mapped onto $L_i\setminus
\{p_j^i\}$ are $A_i\cap B$, hence, using the analysis above and
the fact that every connected component of $B$ must
intersect\footnote{Since $\widetilde{C}_t$ is irreducible, thus
$\widetilde{C}_0$ is connected.} $A_0\cup A_\infty$, we conclude
that $C_0''$ varies in a family of dimension
$$\begin{array}{l}
dim\le-f^*(K_{\Sigma_n}+L_0+L_\infty).B+p_a(B)-1+
\end{array}$$
$$\begin{array}{r}
|A_0\cap B|+|A_\infty\cap B|=
\end{array}$$
$$\begin{array}{l}
-(K_{\Sigma_n}+L_0+L_\infty).C_t+g-2+
\end{array}$$
$$\begin{array}{r} (K_{\Sigma_n}+L_0+L_\infty).(s_0L_0+s_\infty
L_\infty)-p_a(A_0)-p_a(A_\infty)+3.
\end{array}$$
Since $p_a(A_i)\ge 1-s_i$,
$$dim\le-(K_{\Sigma_n}+L_0+L_\infty).C_t+g-2-(s_0+s_\infty)+1.$$

On the other hand $C_0''$ must vary in a family of dimension at
least $-C_t.(K_{\Sigma_n}+L_0+L_\infty)+g-2$, hence $s_0=1$,
$s_\infty=0$, and all the inequalities above are equalities.
Furthermore,
\begin{itemize}
\item Over any irreducible component of $C_0''$ there is a unique
irreducible component of $B$ mapped surjectively onto this
component.
\item $A_0$ is a tree whose root $A_0^R\simeq\PP^1$ mapped
isomorphically onto $L_0$. Any connected component of
$A_0\setminus A_0^R$ intersects at most one connected component of
$B$ and at exactly one point.
\item $C_0''$ is reduced.
\item any two irreducible components of $C_0''$ intersect
transversally.
\item $f(A_0\cap B)$ is a set of generic points of $L_0$, in particular it is disjoint from $p_1^0,...,p_{nd+k}^0$.
\item $C_0''=C_0'$ has geometric genus $g+1-|A_0\cap B|$.
\item $C_0'$ intersects $L_0$ transversally outside of $f(A_0\cap B)$.
\end{itemize}
Now we can describe $C_0''$ explicitly. $C_0''$ is a reduced nodal
curve intersecting $L_\infty$ transversally, and it is smooth at
the points of intersection with $L_0$. And, finally, if $p\in
A_0\cap B$ then in a neighborhood of $f(p)$ the delta invariant of
$C_t$ is equal to the local delta invariant of $C_0$ minus one for
all sufficiently small values of $t$, hence if $m$ denotes the
order of tangency of $C_0'$ and $L_0$ at $f(p)$ then $C_t$ has
$m-1$ nodes in a small neighborhood of $f(p)$. \qed
\end{proof}

To complete the proof we must show that $V$ contains nodal
equige\-neric deformations of $C_0$, since any such deformation
must be of the form $E=L\cup E'$, where $L\equiv L_0$ and
$E'\equiv C'_0$.

Let us denote the points of intersection of $L_0$ with $C_0'$
other than $p_1^0,...,p_{nd+k}^0$ by $q_1,...,q_r$, and the orders
of the tangency by $m_1,...,m_r$ respectively. Thus
$\{q_1,...,q_r\}=f(A_0\cap B)$. Consider now the product
$\CV=\prod_{i=1}^r\CV_i$ of the versal deformations of the
tacnodes of orders $m_1,...,m_r$, and consider the natural map
$$\psi:(V, C_0)\rTo \CV.$$
Then $\psi(V)\subset \prod_{i=1}^r\CV_i(m_i-1),$ where $\CV_i(h)$
denotes the closure of the locus of deformations having $h$ nodes.
\begin{Claim}[\cite{CH98.2} Lemma 2.8]\label{cl:deftac} $\CV_i$ are smooth, irreducible of dimension $2m_i-1$,
and for $m_i-1\le h\le m_i$, $\CV_i(h)$ are irreducible of
dimension $2m_i-1-h$.
\end{Claim}

We denote the nodes of $C_0$ different from $q_1,...,q_r$ by
$o_1,...,o_\delta$. Consider the component $U$ of the decorated
Severi variety containing $(C_0;o_1,...,o_\delta)$. Thus $U$ is
smooth by Proposition \ref{prop:smoothness} (1). Now, let
$$\phi:(U, (C_0;o_1,...,o_\delta))\rTo (\CV, 0)$$ be the natural
map.

\begin{Claim}\label{cl:defsurj}
The map $d\phi :T_{(C_0;o_1,...,o_\delta)}U\rTo T_0\CV$ is
surjective.
\end{Claim}

We postpone the proof of the Claim, and first finish the proof of
Proposition \ref{contain}. Since $U$ and $\CV$ are smooth, Claim
\ref{cl:defsurj} implies that the central fiber $\phi^{-1}(0)$ is
smooth at $(C_0;o_1,...,o_\delta)$, and since the subvariety
$\prod_{i=1}^r\CV_i(m_i-1)\subset \CV$ is irreducible
$$U'=\phi^{-1}\left(\prod_{i=1}^r\CV_i(m_i-1)\right)$$ is also
irreducible. Another conclusion of Claim \ref{cl:defsurj} is the
surjectivity of $\phi$. Thus the generic point of $U'$ corresponds
to a nodal equigeneric deformations of $(C_0;o_1,...,o_\delta)$,
since $\prod_{i=1}^r\CV_i(m_i)\subset\prod_{i=1}^r\CV_i(m_i-1)$.

Consider the natural projection $U'\rTo V_{g,d,k}$. Its image
intersects $V$, and, since $U'$ is irreducible, it belongs to $V$.
Thus $V$ contains nodal equigeneric deformations of $C_0$, and we
are done. \qed

%\begin{corollary}
%$V$ contains nodal curves $E$ of type $E=L\cup E'$, where $L\equiv
%L_0$, and $E'\in V_{g+1-|A_0\cap B|,d,k-1}$.
%\end{corollary}

\begin{proof}[of Claim \ref{cl:defsurj}]
To prove the claim, one must interpret the tangent spaces and the
differential map in cohomological terms. Following the proof of
Proposition \ref{prop:smoothness}, one can see that
$$T_{(C_0;o_1,...,o_\delta)}U\simeq H^0(C_0, \CJ(dL_0+kF)),$$ where
$\CJ$ denotes the ideal sheaf of the zero dimensional scheme
$\cup_{i=1}^\delta o_i\subset C_0$. The tangent space to $\CV$ is
isomorphic to $\oplus_{i=1}^r\CO^{es}(q_i)$, where
$\CO^{es}=\CO_{C_0}/I^{es}$, and $I^{es}$ denotes the
equisingularity ideal of $C_0$. We define $X^{es}=\Spec \CO^{es}$.
Since $A_{2m-1}$ is a simple singularity for any $m\ge 1$ the
equisingularity ideal $I^{es}$ is generated locally at a singular
point by the partial derivatives of the defining equation of the
curve. In these notations the map $d\phi$ is given by the natural
restriction map $H^0(C_0, \CJ(dL_0+kF))\rTo
\oplus_{i=1}^r\CO^{es}(q_i),$ associated to the short exact
sequence
$$0\rTo I^{es}(dL_0+kF)\rTo \CJ(dL_0+kF)\rTo \oplus_{i=1}^r\CO^{es}(q_i)\rTo 0.$$

To prove that $d\phi$ is surjective it is sufficient to show that
\begin{equation}\label{eq:defsurjh1}
H^1(C_0, I^{es}(dL_0+kF))=0.
\end{equation}
It is also a necessary condition since $H^1(C_0, \CJ(dL_0+kF))=0$
by Claim \ref{cl:h1van1}. Consider the short exact sequence of
sheaves
$$0\rTo I_{X^{es}:L_0}\rTo^{\times L_0} I^{es}(L_0)\rTo I_{X^{es}\cap
L_0/L_0}(L_0)\rTo 0.$$ Thus to prove (\ref{eq:defsurjh1}), it is
sufficient to show that
\begin{equation}\label{eq:h101}
H^1(C_0, I_{X^{es}:L_0}((d-1)L_0+kF))=0
\end{equation}
and
\begin{equation}\label{eq:h102}
H^1(L_0, I_{X^{es}\cap L_0/L_0}(dL_0+kF)=0.
\end{equation}
Let $p$ be any point in the support of $X^{es}$. If $p\notin L_0$
then $$I_{X^{es}:L_0, p}=I^{es}_p=\CJ^{cond}_p\;\; \mbox{and}\;\;
I_{X^{es}\cap L_0/L_0, p}=0,$$ where $\CJ^{cond}$ denotes the
conductor ideal of $C_0$. If $p\in L_0$ then consider a local
system of coordinates $x,y$ at $p\in \Sigma_n$ such that $L_0$ is
given by $y=0$, and $C'_0$ is given by $y=x^m$. In these notations
$$I_{X^{es}:L_0, p}=y\CO_{C_0,p}+x^{m-1}\CO_{C_0,p}\supset y\CO_{C_0,p}+x^m\CO_{C_0,p}=\CJ^{cond}_p$$
and $$I_{X^{es}\cap L_0/L_0, p}=x^m\CO_{L_0,p}.$$ Thus
$I_{X^{es}\cap L_0/L_0}(dL_0+kF)$ is a line bundle of degree
$C_0.L_0-C'_0.L_0=n\ge 0$, hence $H^1(L_0, I_{X^{es}\cap
L_0/L_0}(dL_0+kF)=0$ by the Riemann-Roch theorem, which proves
(\ref{eq:h102}).

It follows from the description above that $\CJ^{cond}\subset
I_{X^{es}:L_0}$. Thus $H^1(C_0, \CJ^{cond}((d-1)L_0+kF))=0$
implies (\ref{eq:h101}). Consider the normalization
$\nu:\widetilde{C}_0\rTo C_0$. It is sufficient to show that
\begin{equation}\label{eq:defsurjh2}
H^1(C, \CJ^{cond}((d-1)L_0+kF)_{|_C})=0,
\end{equation}
for any irreducible component $C\subset \widetilde{C}_0$. Since
$\nu(C)\ne L_\infty$ by Lemma \ref{lm:rednotcont}, thus
$(2F+L_\infty).\nu(C)>0$ and $\CJ^{cond}((d-1)L_0+kF)$ restricted
to $C$ is a line bundle of degree
$$\begin{array}{ll}\nu(C).C'_0-\deg(\CJ^{cond}_{|_C})&=\nu(C)^2-\nu(C).L_0-2\delta(\nu(C))\\&=2g(C)-2-(K_{\Sigma_n}+L_0).\nu(C)\\&=2g(C)-2+(2F+L_\infty).\nu(C)>2g(C)-2,\end{array}$$
where $\delta(\nu(C))$ denotes the number of nodes of $\nu(C)$
(which is equal to the total delta invariant of $\nu(C)$), thus
(\ref{eq:defsurjh2}) follows from the Riemann-Roch theorem. \qed
\end{proof}

\subsection{Proof of Proposition \ref{unique}}
We start with some combinatorics.
\begin{definition}(1) An ordered subset $\mu\subseteq \Gamma^{sing}$ consisting of $r$ nodes is called an $r-$marking
on the curve $\Gamma$.

(2) An $r-$marking $\mu$ is called irreducible if and only if the curve $\Gamma\backslash\mu$ is connected.

(3) We define D-moves on the set of $r-$markings as follows: let
$D, D'\subset\Gamma$ be two different irreducible components, let
$q,q'\in D\cap D'$ be two nodes, and let $\mu=\{p_1,...,p_r\}$ be
an $r-$marking. Then
$$D_{q,q'}(\mu)=\begin{cases}
    \{p_1,...,p_{i-1},q,p_{i+1},...,p_r\} & \text{if $q\notin\mu$ and $q'=p_i$}, \\
    \{p_1,...,p_{i-1},q',p_{i+1},...,p_r\} & \text{if $q=p_i$ and $q'\notin\mu$}, \\
    \{p_{\tau_{ij}(1)},...,p_{\tau_{ij}(r)}\} & \text{if $q=p_i$ and $q'=p_j$}, \\
    \mu & \text{otherwise,}
  \end{cases}
$$
where $\tau_{ij}\in {\mathfrak S}_r$ denotes the elementary transposition $\tau_{ij}=(i\,j)$.

(4) Assume that $n>0$. We define T-moves on the set of
$r-$mar\-kings as follows: let $D, D', D''\subset\Gamma$ be three
different irreducible components, and let $q\in D'\cap D''$,
$q'\in D\cap D''$, $q''\in D\cap D'$ be three nodes, and let
$\mu=\{p_1,...,p_r\}$ be an $r-$marking. Then if $q'\notin\mu$ we
define
$$T_{q,q',q''}(\mu)=
  \begin{cases}
    \{p_1,...,p_{i-1},q'',p_{i+1},...,p_r\} & \text{if $q=p_i; q''\notin\mu$}, \\
    \{p_1,...,p_{i-1},q,p_{i+1},...,p_r\} & \text{if $q''=p_i; q\notin\mu$}, \\
    \{p_{\tau_{ij}(1)},...,p_{\tau_{ij}(r)}\} & \text{if $q=p_i, q''=p_j$}, \\
    \mu & \text{if $q,q''\notin\mu$,}
  \end{cases}
$$
otherwise we define $T_{q,q',q''}(\mu)=\mu$.

(5) Assume that $n=0$. We define Q$^h$-moves and Q$^v$-moves as
follows: let $\mu=\{p_1,...,p_r\}$ be an $r-$marking, let $X, X',
Y, Y'\subset\Gamma$ be four different irreducible components
satisfying $X\equiv X'\equiv L_0$ and $Y\equiv Y'\equiv F$, and
let $q\in X\cap Y$, $q'\in X\cap Y'$, $q''\in X'\cap Y$, $q'''\in
X'\cap Y'$ be nodes. If $q'',q'''\notin\mu$ we define
$$Q^h_{q,q',q'',q'''}(\mu)=
  \begin{cases}
    \{p_1,...,p_{i-1},q',p_{i+1},...,p_r\} & \text{if $q=p_i$ and $q'\notin\mu$}, \\
    \{p_1,...,p_{i-1},q,p_{i+1},...,p_r\} & \text{if $q'=p_i$ and $q\notin\mu$}, \\
    \{p_{\tau_{ij}(1)},...,p_{\tau_{ij}(r)}\} & \text{if $q=p_i$ and $q'=p_j$}, \\
  \end{cases}
$$
otherwise, we define $Q^h_{q,q',q'',q'''}(\mu)=\mu$. If
$q',q'''\notin\mu$, we define
$$Q^v_{q,q',q'',q'''}(\mu)=
  \begin{cases}
    \{p_1,..,p_{i-1},q'',p_{i+1},..,p_r\} & \text{if $q=p_i$ and $q''\notin\mu$}, \\
    \{p_1,..,p_{i-1},q,p_{i+1},..,p_r\} & \text{if $q''=p_i$ and $q\notin\mu$}, \\
    \{p_{\tau_{ij}(1)},..,p_{\tau_{ij}(r)}\} & \text{if $q=p_i$ and $q''=p_j$}, \\
  \end{cases}
$$
otherwise, we define $Q^v_{q,q',q'',q'''}(\mu)=\mu$.

(6) Two $r-$markings, $\mu$ and $\mu'$, are called equivalent if
and only if one can be obtained from another by a sequence of
T-moves, Q-moves, and D-moves.
\end{definition}

\begin{notation} Let $\mu=\{p_1,...,p_r\}$ be any $r-$marking and let $C,C'\subset\Gamma$ be two different components. The following notation will be
useful:
$$\mu_{C,C'}=|\mu\cap C\cap C'|,\:\:\mu_C=|\mu\cap C|,\:\:\mbox{and}\:\: \mu_i=p_i.$$
\end{notation}
\begin{Claim}\label{cl:comb}
Let $r>0$ be such an integer, that the set of irreducible
$r-$markings on the curve $\Gamma$ is not empty. Then for any pair
of distinct irreducible components $C,C'\subset\Gamma$ and for any
$q\in C\cap C'$, there exist irreducible $r-$markings $\mu$ and
$\mu'$ such that $q\in\mu$ and $q\notin\mu'$.
\end{Claim}
\begin{proof} Obvious.
\qed
\end{proof}

From now on we will assume that $n>0$. The remaining case, $n=0$,
is much easier, and the proof in this case can be obtained via the
same lines as in the case we consider. Thus, we leave it to the
reader.

\begin{lemma}\label{l:lcomb}
If $n>0$ then any two irreducible $r-$markings on the curve
$\Gamma$ are equivalent.
\end{lemma}
\begin{proof} It is enough to prove the lemma for the case
$r=k(d-1)+n\frac{d(d-1)}{2}-(d-1)$, since any irreducible
$r-$marking $\mu$ can be extended to an irreducible marking
$\mu^{ex}$ of order $k(d-1)+n\frac{d(d-1)}{2}-(d-1)$, and for any
$D$ or $T$ move $M$ on the extended marking the natural forgetful
map takes $M(\mu^{ex})$ to $M(\mu)$.

Let us prove the lemma by induction on $d+k$. If either $d=2=k+2$
or $d\le 1$, then the lemma is obvious. Assume that the statement
is true for all $d+k\le m$ and let us prove it for $d+k=m+1$. We
can assume that $m>2$ and $d\ge 2$. Let $\mu$ be an irreducible
$r-$marking.

{\it Step 1:} The goal of this step is to prove that there exists a marking $\mu'\sim\mu$ such
that
\begin{equation}\label{combprop}
\mu'_{C,C'}=
  \begin{cases}
    C.C'-1 & \text{if $L_d\in\{C,C'\}$}, \\
    C.C' & \text{otherwise}.
  \end{cases}
\end{equation}
Choose a component $D\subset\Gamma$ in the following way: if $k\ne
0$ then $D=F_k$, otherwise $D=L_1$. Then any irreducible component
$D'\subset\Gamma$ different from $D$ and satisfying $D.D'>0$
belongs to the linear system $|\CO_{\Sigma_n}(L_0)|$. Moreover,
there are at least two such components. Now, let us choose
$\widetilde{\mu}\sim\mu$, such that
$$\widetilde{\mu}_D=\max_{\mu'\sim\mu}\{\mu'_D\}.$$
Then $\widetilde{\mu}_{|_{\overline{\Gamma\backslash D}}}$ is an
irreducible marking on $\overline{\Gamma\backslash D}$, since
otherwise, due to the choice of $D$, we would be able to find two
distinct irreducible components $D',D''\subset \Gamma$ different
from $D$ such that $\widetilde{\mu}_{D,D'}<D.D'$,
$\widetilde{\mu}_{D,D''}<D.D''$, and
$\widetilde{\mu}_{D',D''}=D'.D''>0$. Hence, there would exist
$q''\in D\cap D'$, $q'\in D\cap D''$, and $q\in D'\cap D''$, such
that $q',q''\notin\mu$ and $q=\widetilde{\mu}_i$ for some $i$.
Thus $\widetilde{\mu}_D<\mu'_D$, where
$\mu'=T_{q,q',q''}(\widetilde{\mu})\sim\mu$, which would be a
contradiction.

Next, we shall prove that
\begin{equation}\label{eqcomblem}
\widetilde{\mu}_D=D.(\Gamma-D)-1.
\end{equation}
Consider two cases: $k>0$ and $k=0$.

\noindent{\it Case 1: $k>0$.} In this case $D=F_k$, and
$$\widetilde{\mu}_D\ge r-\left((k-1)(d-1)+n\frac{d(d-1)}{2}-(d-1)\right),$$
since $\widetilde{\mu}_{|_{\overline{\Gamma\backslash D}}}$ is an
irreducible marking on $\overline{\Gamma\backslash D}$. However,
$$r-\left((k-1)(d-1)+n\frac{d(d-1)}{2}-(d-1)\right)=d-1=D.(\Gamma-D)-1,$$
and we are done.

\noindent{\it Case 2: $k=0$.} In this case $D=L_1$, and
$$\widetilde{\mu}_D\ge r-\left(n\frac{(d-1)(d-2)}{2}-(d-2)\right),$$
since $\widetilde{\mu}_{|_{\overline{\Gamma\backslash D}}}$ is
irreducible. However,
$$r-\left(n\frac{(d-1)(d-2)}{2}-(d-2)\right)=n(d-1)-1=D.(\Gamma-D)-1.$$
Thus, $\widetilde{\mu}_D=D.(\Gamma-D)-1$ in both cases.

Now we can complete the proof of the first step. By
(\ref{eqcomblem}) there exists a unique component
$D_1\subset\Gamma$ satisfying $\widetilde{\mu}_{D,D_1}=D.D_1-1$,
and $\widetilde{\mu}_{D,D'}=D.D'$ for all $D'\subset\Gamma$
different from $D_1$. By the induction hypothesis it remains to
prove that there exists $\mu'\sim\widetilde{\mu}$ such that
$\mu'_{D,L_d}=D.L_d-1$ and $\mu'_D=\widetilde{\mu}_D$. If
$D_1=L_d$ then there is nothing to prove. Otherwise, by the
induction assumption and Claim \ref{cl:comb}, we can find an
irreducible marking $\widehat{\mu}\sim \widetilde{\mu}$ with the
following properties:
$\widehat{\mu}_{D,D'}=\widetilde{\mu}_{D,D'}$ for all
$D'\subset\Gamma$, and $\widehat{\mu}_{D_1,L_d}<D_1.L_d$. Then
there exist $q''\in D\cap D_1$, $q\in D\cap L_d$, and $q'\in
D_1\cap L_d$ such that $q',q''\notin \widehat{\mu}$, and
$q=\widehat{\mu}_i$ for some $i$. Thus
$\mu'=T_{q,q',q''}(\widehat{\mu})$ satisfies the required
condition.

{\it Step 2:} The goal of this step is to prove that any two
markings satisfying (\ref{combprop}) are equivalent. Let $\mu,
\mu'$ be two such markings. Applying several D-moves we can find a
marking $\mu''\sim \mu'$ such that $\mu''$ differs from $\mu$ only
by the order of marked points. Namely, there exists $\tau\in
{\mathfrak S}_r$ such that $\mu''_i=\mu_{\tau(i)}$ for all $i$. We
use the notation $\tau(\mu)$ for such $\mu''$. It remains to prove
that $\mu\sim \tau(\mu)$. Without loss of generality we can assume
that $\tau=\tau_{ij}$ is a simple transposition; moreover we can
assume that $i=1$ and $j=2$. Let $D_1,D_2,D_3,D_4$ be components
such that $\mu_1\in D_1\cap D_2$ and $\mu_2\in D_3\cap D_4$. If
$\{D_1, D_2\}=\{D_3, D_4\}$, then we apply  $D_{\mu_1,\mu_2}$ to
finish the proof. So we can assume that $D_1\ne D_3, D_4$ and
$D_2\ne D_3$.

If $\mu_1\in L_d\cap D'$ for some $D'$, then $D'\equiv L_0$ due to
(\ref{combprop}). Thus there exists a component $D''\ne D', L_d$
such that $D''.D'>0$. Let $\mu_i\in D''\cap D'$ be any node, and
let $q\in L_d\cap D''$ be the node not belonging to $\mu$. Then
$\mu\sim T_{\mu_1,q,\mu_i}(\mu)$, and hence we can reduce the
statement to the case when $\mu_1\notin L_d$. Applying the same
argument to $\mu_2$ we get one of the following: either
$\mu_i=\mu_1$, and then $\mu\sim T_{\mu_2,q,\mu_i}(\mu)=\mu''$, or
$\mu_1,\mu_2\notin L_d$, and thus that $D_i\ne L_d$ for all $i$.
Now we shall consider several cases:

\noindent {\it Case 1: $D_2=D_4\equiv L_0$.}
 Let $q_1\in L_d\cap D_1$, $q_2\in L_d\cap D_2$, and $q_3\in L_d\cap D_3$ be the nodes not belonging
to $\mu$. Then $$\tau_{12}(\mu)=T_{\mu_2,q_3,q_2}(T_{\mu_1,q_1,q_2}(T_{\mu_2,q_3,q_2}(\mu)))\sim\mu.$$

\noindent {\it Case 2: $D_2=D_4\equiv F$.} In this case $D_1\equiv
D_3\equiv L_0$. Let $\mu_i\in D_1\cap D_3$ be any node. Then
$\tau_{12}=\tau_{1i}\circ\tau_{i2}\circ\tau_{1i}$, and thus we
reduce to the previous case. So $\tau_{12}(\mu)\sim\mu$.

\noindent {\it Case 3: $\{D_1, D_2\}\cap\{D_3, D_4\}=\emptyset$.}
Without loss of generality $D_1\equiv D_3\equiv L_0$. Let
$\mu_i\in D_1\cap D_3$ be any node. Then
$\tau_{12}=\tau_{1i}\circ\tau_{i2}\circ\tau_{1i}$, so
$\tau_{12}(\mu)\sim\mu,$ by the first case, and we are done, since
cases 1, 2, 3 cover all the possibilities. \qed
\end{proof}

\begin{lemma}\label{cl:irrp1}
Let $X=X_1\cup X_2\cup...\cup X_r$ be a nodal curve. Assume that
$X_2,...,X_r$ are generic curves of types $L_0$ and $F$, and
assume that $X_1$ is a generic nodal rational curve whose type is
one of $L_0$, $L_0+F$ or $2L_0$.
%For any $\CL\in Pic(\Sigma_n)$ define $s=\deg(\CL\otimes\CO_X)$, and
Consider the variety
$W(X, \CL)\subset (X^{smooth})^s\times |\CL|$ given by
$$\left\{(p_1,...,p_s;\xi)\,:\,p_i\ne p_j\, \mbox{for all}\:\: i\ne j,\,\mbox{and}\:\: \xi(p_1)=...=\xi(p_s)=0\right\},$$
where $s=\deg(\CL\otimes\CO_X)$ and either $\CL\simeq
\CO_{\Sigma_n}(L_0)$ or $\CL\simeq \CO_{\Sigma_n}(F)$. Let
$(p_1,...,p_s;\xi)\in W(X, \CL)$ be an arbitrary point, and let
$W\subseteq W(X, \CL)$ be the irreducible component containing
$(p_1,...,p_s;\xi)$. Assume that $n>0$ and there exist $a\neq b$
satisfying $p_a,p_b\in X_1$. Then
$$(p_{\tau_{ab}(1)},...,p_{\tau_{ab}(s)};\xi)\in W,$$
where $\tau_{ab}\in {\mathfrak S}_s$ denotes the elementary transposition $\tau_{ab}=(a\,b)$.
\end{lemma}
\begin{proof} If $\CL\simeq \CO_{\Sigma_n}(F)$, we can assume that $X_1\equiv 2L_0$.
Consider the variety $W(X_1, \CL)$, which is irreducible since
$X_1$ is irreducible and the projection to the first factor
$W(X_1,\CL)\rTo X_1$ is dominant and has irreducible fibers.
%and $s=2$.
The natural forgetful map
$f:W\rTo W(X_1, \CL)$ is dominant, and since $F.X_i\le 1$ for all
$i>1$ it is also one-to-one. This implies the statement.

Assume now that $\CL\simeq \CO_{\Sigma_n}(L_0)$. It is easy to see
that $W(X, \CL)$ is smooth. Let $(x_1,...,x_s;\eta)\in
\overline{W(X, \CL)}$ be any point with the following properties:
$x_i\in X^{smooth}$ for all $i$, $x_a=x_b$,
$\eta(x_a)=d(\eta_{|_X})(x_a)=0$, and $x_i\ne x_j$ for all $\{i,
j\}\ne \{a, b\}$. Then $(x_1,...,x_s;\eta)$ is a smooth point of
$\overline{W(X, \CL)}$. Thus to prove the lemma it suffices to
show that $\overline{W}$ contains a point with such properties.
Consider the forgetful map $f:\overline{W}\rTo (X_1^{smooth})^2$
given by $f(x_1,...,x_s;\eta)=(x_a,x_b)$. This map is surjective,
since $\CL\simeq \CO_{\Sigma_n}(L_0)$ and $n>0$, and it is clear
that a generic $\alpha\in f^{-1}(\Delta_{X_1})$ satisfies the
properties mentioned above. \qed
\end{proof}

\begin{proof}[of Proposition \ref{unique}] We start with
the following remark: let $U\subset U_{d,k,\delta}$ be an
irreducible component containing $(\Gamma; p_1,...,p_\delta)$. We
define a $\delta-$marking $\mu^{\Gamma;p_1,...,p_\delta}$ on the
curve $\Gamma$ in the following way:
$$\mu^{\Gamma;p_1,...,p_\delta}=\{p_1,...,p_\delta\}.$$ Then the generic curve $C\in
\overline{\psi(U)}$ is irreducible if and only if the marking
$\mu^{\Gamma;p_1,...,p_\delta}$ is irreducible. The collection of
all $\delta-$markings corresponding to $U$ is denoted by
$$\CM(U)=\{\mu^{\Gamma;p_1,...,p_\delta}\,|\,(\Gamma;
p_1,...,p_\delta)\in U\}.$$ Thus Proposition would be proven once
we show that $\CM(U)$ is closed under T-moves and under D-moves.
This is what we proceed to proving now under the assumption $n>0$.

{\it Step 1:} First, we shall prove that $\CM(U)$ is closed under
T-moves. We choose an arbitrary marking
$\mu=\mu^{\Gamma;p_1,...,p_\delta}\in \CM(U)$ and label the rest
of the nodes of $\Gamma$ by $p_{\delta+1},...,p_{\delta'}$, where
$\delta'=dk+n\frac{d(d-1)}{2}$. Let $D,D',D''\subset\Gamma$ be
three different irreducible components, and let $q\in D'\cap D''$,
$q'\in D\cap D''$, $q''\in D\cap D'$ be three nodes. If $q'\in\mu$
then $T_{q,q',q''}(\mu)=\mu\in\CM(U)$ and we are done. So we can
assume that $q'\notin\mu$. If $q,q''\notin\mu$ then again
$T_{q,q',q''}(\mu)=\mu\in\CM(U)$ and we are done. So without loss
of generality we can assume that $q\in\mu$. We shall show that
$T_{q,q',q''}(\mu)\in\CM(U)$. Without loss of generality, $q=p_1$,
$q''=p_i$, and $q'=p_{\delta'}$.

Consider the irreducible component $U'\subset U_{d,k,\delta'-1}$
containing the pointed curve $(\Gamma; p_1,...,p_{\delta'-1})$,
and let $f:U'\rTo U$ be the natural forgetful map. It is
sufficient to prove that
$$(\Gamma; p_i,p_2,...,p_{i-1},p_1,p_{i+1},...,p_{\delta'-1})\in U',$$
since $T_{q,q',q''}(\mu)=\mu^{f(\Gamma;
p_i,p_2,...,p_{i-1},p_1,p_{i+1},...,p_{\delta'-1})}$.

Let $(C;x_1,...,x_{\delta'-1})\in U'$ be a generic element. Then
$C$ has a unique component $C_2$ of type $D+D''$ among its $d+k-1$
irreducible components. Moreover, there exists another irreducible
component $C_1$ such that $x_1, x_i\in C_1\cap C_2$. We denote
$C^{fix}=\cup_{l=2}^{d+k-1}C_l$, where $C_1,...,C_{d+k-1}$ are the
irreducible components of $C$.

Consider the locus $U''\subset U'$ of pointed curves
$(C';x'_1,...,x'_{\delta'-1})$ with the following property:
$C'=C'_1\cup C^{fix}$, where $C'_1\equiv C_1$ is generic. Let
$U'''\subset U''$ be the irreducible component containing the
pointed curve $(C;x_1,...,x_{\delta'-1})$, and let
$$\phi:U'''\rTo W(C^{fix},\CO_{C^{fix}}(C_1))$$
be the natural forgetful map (cf. Lemma \ref{cl:irrp1}). Then the
image of $\phi$ is dominant in an irreducible component of
$W(C^{fix},\CO_{C^{fix}}(C_1))$, and $\phi$ is one-to-one. Hence
$$(\Gamma; p_i,p_2,...,p_{i-1},p_1,p_{i+1},...,p_{\delta'-1})\in
U'''\subset U'$$ by Lemma \ref{cl:irrp1}.

{\it Step 2:} The goal of this step is to prove that $\CM(U)$ is
closed under D-moves. We choose an arbitrary marking
$\mu=\mu^{\Gamma;p_1,...,p_\delta}\in \CM(U)$ and label the rest
of the nodes of $\Gamma$ by $p_{\delta+1},...,p_{\delta'}$, where
$\delta'=dk+n\frac{d(d-1)}{2}$. Let $D,D'\subset\Gamma$ be two
different irreducible components, and let $q,q'\in D\cap D'$ be
two nodes. If $q,q'\notin\mu$ then $D_{q,q'}(\mu)=\mu\in\CM(U)$
and we are done. Thus, without loss of generality, we can assume
that $q\in\mu$. We shall show that $D_{q,q'}(\mu)\in\CM(U)$.

Consider the irreducible component $U'\subset U_{d,k,\delta'}$
containing the pointed curve $(\Gamma; p_1,...,p_{\delta'})$, and
let $f:U'\rTo U$ be the natural forgetful map. Then
$$D_{q,q'}(\mu)=\mu^{f(\Gamma; p_{\tau_{ij}(1)},...,p_{\tau_{ij}(\delta')})},$$
where $q=p_i$, $q'=p_j$, and $\tau_{ij}\in {\mathfrak
S}_{\delta'}$ denotes the elementary transposition
$\tau_{ij}=(i\,j)$. Hence it is sufficient to prove that
$$(\Gamma; p_{\tau_{ij}(1)},...,p_{\tau_{ij}(\delta')})\in U'.$$
Let $\Gamma_1,...,\Gamma_{d+k}$ be the irreducible components of
$\Gamma$. We can assume that $D=\Gamma_1$, and we denote
$\Gamma^{fix}=\cup_{l=2}^{d+k}\Gamma_l$.

Consider the locus $U''\subset U'$ of pointed curves
$(C;x_1,...,x_{\delta'})$ with the following property: $C=C_1\cup
\Gamma^{fix}$, where $C_1\equiv \Gamma_1$ is generic. Let
$U'''\subset U''$ be the irreducible component containing
$(\Gamma;p_1,...,p_{\delta'})$, and let
$$\phi:U'''\rTo W(\Gamma^{fix},\CO_{\Gamma^{fix}}(\Gamma_1))$$
be the natural forgetful map (cf. Lemma \ref{cl:irrp1}). Then the
image of $\phi$ is dominant in an irreducible component of
$W(\Gamma^{fix},\CO_{\Gamma^{fix}}(\Gamma_1))$, and $\phi$ is
one-to-one. Hence
$$(\Gamma; p_{\tau_{ij}(1)},...,p_{\tau_{ij}(\delta')})\in
U'''\subset U'$$ by Lemma \ref{cl:irrp1}. \qed
\end{proof}

\section{Rational curves on toric surfaces}

The goal of this section is to prove Conjecture \ref{conj0} for
the case of rational curves.

\begin{proposition} Let $\Sigma=Tor(\Delta)$ be a toric surface assigned to an integral polygon $\Delta\subset \RR^2$,
and let $\CL\in Pic(\Sigma)$ be an effective class. Consider
variety $V$ parameterizing all irreducible nodal rational curves
in the linear system $|\CL|$ belonging to the smooth locus of
$\Sigma$. Then $V$ is either empty or irreducible.
\end{proposition}

\begin{proof}
We can resolve the singularities of $\Sigma$ by a sequence of blow
ups of the singular zero-dimensional orbits. Since $V$
parameterizes curves that do not contain singularities of
$\Sigma$, the variety $V$ parameterizes also irreducible rational
curves in the pull back of $\CL$ to the disingularization of
$\Sigma$. Thus to prove the Proposition, it is sufficient to
consider only the case of smooth surface $\Sigma$. So let us
assume that $\Sigma$ is smooth.

Let $(a_i,b_i)$, $1\le i\le n$, be the primitive integral vectors
parallel to the sides of the $n$-gon $\Delta$ oriented
counterclockwise. We define $(a_{n+1},b_{n+1})=(a_1,b_1)$. Then
$\{(a_i,b_i),(a_{i+1},b_{i+1})\}$ is a basis of the integral
lattice for any $1\le i\le n$, since $\Sigma$ is smooth.

Now let $C\in V$ be a generic element, and let $\phi:\PP^1\to
\Sigma$ be a parameterization of $C$. If $C$ coincides with one of
the boundary components then $V$ is a point and we are done. Thus
we can assume that $C$ intersects the boundary divisor at a finite
number of points, moreover there are at least two such points,
since no chart isomorphic to $\KK^2$ can contain a complete curve.
Thus the first chern class of the normal bundle $\CN_\phi$ is
non-negative, and, since $C$ is nodal, $\CN_\phi$ is a line
bundle. So, we can conclude that first $V$ is equidimensional,
and, second, no irreducible component of $V$ has a fixed point, in
particular $C$ contains no zero-dimensional orbits.

For any $1\le i\le n$ we define $\{c_{ij}\}=\phi^{-1}(L_i)\subset
\PP^1$, where $L_i\subset \Sigma$ is the one-dimensional orbit
corresponding to $(a_i,b_i)$. We can assume that $\phi(\infty)\in
(\KK^*)^2$. The restriction of $\phi$ to ${\mathbb
A}^1\setminus\{c_{ij}\}$ is given by two invertible functions
$x(t),y(t)\in\KK[t,(t-c_{ij})^{-1}]$, hence
$$x(t)=\alpha\prod (t-c_{ij})^{m_{ij}}, \: \mbox{and}\: \: y(t)=\beta\prod (t-c_{ij})^{n_{ij}}.$$

Let $1\le i\le n$ be any index such that $k_i=|\phi^{-1}(L_i)|>0$.
Consider the affine plane
$\Spec\KK[x^{-a_i}y^{-b_i},x^{a_{i+1}}y^{b_{i+1}}]\subset \Sigma$.
In this chart the line $L_i$ is given by
$x^{a_{i+1}}y^{b_{i+1}}=0$ and $x^{-a_i}y^{-b_i}\ne 0$, hence
$a_{i+1}m_{ij}+b_{i+1}n_{ij}=k_{ij}$ and $a_im_{ij}+b_in_{ij}=0$
for all $j$, where $k_{ij}>0$ denotes the order of $\phi^*(L_i)$
at $c_{ij}$. Since $\{(a_i,b_i),(a_{i+1},b_{i+1})\}$ is a
(positive) basis of the integral lattice, we can conclude that
$n_{ij}=k_{ij}a_i$ and $m_{ij}=-k_{ij}b_i$ for all $i$ and $j$.
Thus
$$x(t)=\alpha\prod (t-c_{ij})^{-k_{ij}b_i},\: \mbox{and}\: \: y(t)=\beta\prod (t-c_{ij})^{k_{ij}a_i}.$$
We shall mention that $\deg x(t)=\deg y(t)=0$ since
$\phi(\infty)\in (\KK^*)^2$.

Next we would like to show that $k_{ij}=1$ for all $i$ and $j$. If
this is not the case, then without loss of generality we can
assume that $k_{11}>1$. Thus the locus of rational curves (in the
same linear system) that admit the following parameterization
$$x(t)=\alpha(t-c'_{11})^{-(k_{11}-1)b_1}(t-c''_{11})^{-b_1}\prod_{(i,j)\ne(1,1)}(t-c_{ij})^{-k_{ij}b_i},$$
$$y(t)=\beta(t-c'_{11})^{(k_{11}-1)a_1}(t-c''_{11})^{a_1}\prod_{(i,j)\ne(1,1)}(t-c_{ij})^{k_{ij}a_i},$$
has dimension grater than $\dim V$, which is a contradiction.

Now we see that $V$ contains an open dense subset isomorphic to an
open subset of
$$(\KK^*)^2\times\prod_{i=1}^n Sym^{k_i}\PP^1$$ modulo the automorphisms of
$\PP^1$. Thus $V$ contains an irreducible dense open subset, hence
$V$ is irreducible. \qed
\end{proof}

Finally we would like to explain why the assumption that a generic
$C\in V$ does not contain singularities of $\Sigma$, is necessary.
Consider the toric surface $\Sigma$ assigned to the triangle
$\{(0,0),(0,2),(4,0)\}$, and let $\CL$ be the tautological line
bundle on $\Sigma$. This surface has unique singular point. Then
the locus of irreducible rational curves in $|\CL|$ consists of
two irreducible components of dimension 7. The first component was
described in the proposition. To see the second component let us
consider the desingularization of $\Sigma$, which is isomorphic to
the Hirzebruch surface $\Sigma_2$ assigned to the trapezia
$\{(0,0),(0,1),(2,1),(4,0)\}$. Then the tautological linear system
on $\Sigma_2$ has dimension 7 and the projection of this system to
$\Sigma$ defines a hyperplane in $|\CL|$ consisting of curves
passing through the node of $\Sigma$. This hyperplane is the
second component of the Severi variety.

\end{document}